\pdfoutput=1
\newcommand{\online}{0}
\newcommand{\classstyle}{1}
\documentclass[12pt]{article}
\usepackage[utf8]{inputenc} 
\usepackage{fullpage}
\usepackage{times}
\usepackage[normalem]{ulem}
\usepackage{fancyhdr,graphicx,amsmath,amssymb, mathtools, scrextend, titlesec, enumitem}
\usepackage{bm}
\usepackage{tcolorbox}
\usepackage{amsthm}
\usepackage{mathrsfs} 
\usepackage[acronym,nomain,nonumberlist]{glossaries}
\usepackage{siunitx}
\makeglossaries

\usepackage{bbding}

\usepackage{caption}
\usepackage{subcaption}

\usepackage{tikz}

\usepackage{amsfonts}
\numberwithin{equation}{section}
\ifnum\classstyle=3
\else

\fi

\ifnum\classstyle=3
\else
\usepackage[hypertexnames=false]{hyperref}
\hypersetup{
    colorlinks = true,
    linkbordercolor = {white},
    linkcolor = {blue!50!black},
citecolor     = {blue!50!black},
}
\fi

\usepackage{xcolor,colortbl}

\usepackage{nicematrix}

\usepackage{algorithm}
\usepackage{algpseudocode}

\algrenewcommand\algorithmicindent{2.0em}

\usepackage{tikz}
\usepackage{pdfpages}

\usepackage{chngcntr}

\counterwithin{figure}{section}

\usepackage{ulem}

\usepackage{yhmath}

\usepackage{accents}
\newlength{\dhatheight}

\makeatletter
\@namedef{subjclassname@2020}{%
  \textup{2020} Mathematics Subject Classification}
\makeatother

\usepackage{arxiv}
\usepackage[style=ieee,sorting=none,natbib=true,maxbibnames=99]{biblatex}
\DeclareNameAlias{sortname}{family-given}
\DeclareNameAlias{default}{family-given}
\usepackage{cleveref}
\usepackage{setspace}
\usepackage{indentfirst}
\title{A two-step surrogate method for sequential uncertainty quantification in high-dimensional inverse problems}

\author{
 \textbf{Ningxin~Yang}\thanks{corresponding author: n.yang23@imperial.ac.uk}\\
Department of Civil \\ and Environmental Engineering \\ Imperial College London\\  London, UK\\
 	\And
	\textbf{Truong~Le}\\
Department of Civil \\ and Environmental Engineering \\ Imperial College London\\  London, UK\\ 
 \And
 \textbf{Lidija~Zdravković} \\
	Department of Civil \\ and Environmental Engineering \\ Imperial College London\\  London, UK\\
	\And
 \textbf{David~Potts} \\
Department of Civil \\ and Environmental Engineering \\ Imperial College London\\  London, UK
}

\renewcommand{\shorttitle}{}

\begin{document}
\maketitle

\begin{abstract}
Predictive estimation, which comprises model calibration, model prediction, and validation, is a common objective when performing inverse uncertainty quantification (UQ) in diverse scientific applications. These techniques typically require thousands to millions of realisations of the forward model, leading to high computational costs. Surrogate models are often used to approximate these simulations. However, many surrogate models suffer from the fundamental limitation of being unable to estimate plausible high-dimensional outputs, inevitably compromising their use in the UQ framework. To address this challenge, this study introduces an efficient surrogate modelling workflow tailored for high-dimensional outputs. Specifically, a two-step approach is developed: (1) a dimensionality reduction technique is used for extracting data features and mapping the original output space into a reduced space; and (2) a multivariate surrogate model is constructed directly on the reduced space. The combined approach is shown to improve the accuracy of the surrogate model while retaining the computational efficiency required for UQ inversion. The proposed surrogate method, combined with Bayesian inference, is evaluated for a civil engineering application by performing inverse analyses on a laterally loaded pile problem. The results demonstrate the superiority of the proposed framework over traditional surrogate methods in dealing with high-dimensional outputs for sequential inversion analysis. 


\end{abstract}

\section{Introduction}\noindent
An inverse problem is where the input, or part of it, to a problem is sought, knowing some of the results (outputs). For example, in civil engineering, one might be faced with the construction of an excavation which has restrictions on movements due to adjacent services and infrastructure. In the design stage, forward calculations are performed using advanced numerical analysis to predict movements. However, these will be based on assumed material properties, initial ground conditions, and approximate boundary conditions, which will involve some uncertainties \cite{jung2009,jin2021}. Once the excavation commences, field measurements of movements could become available and are unlikely to match initial predictions. Consequently, an inverse problem arises, requiring the use of these measurements to refine estimates of material properties and initial conditions. These refined estimates can then be used in forward calculations to update movement predictions for later excavation stages. Such calculations can also be used to quantify uncertainty and to establish the most influential soil parameters, initial conditions and/or boundary conditions. 

In many engineering problems, inverse analysis involves performing many forward calculations to determine the best set of model parameters. Given the numerical nature and the diversity of parameters involved, Monte Carlo sampling is necessary to form the set of forward calculations. 
The number of calculations required increases with the number of input parameters (material properties and initial boundary conditions) and the number of output quantities, often necessitating thousands to millions of forward calculations. In many cases, especially where the output quantity consists of a field variable (e.g. stress or displacement field), the computational cost of performing such a large number of forward calculations becomes prohibitive for real world practices. To overcome this challenge and facilitate progress, a surrogate model can be established. This model mimics the behavior of the original numerical analysis tool but performs forward calculations at a significantly reduced computational cost.

To reduce the computational cost further, it is often possible to identify input conditions (material properties and/or initial and boundary conditions) which have a dominant effect on the outputs of interest. Likewise, certain output quantities (e.g., movement at a certain location) may have a dominant effect on all other output quantities of interest (i.e., there is dependency between the outputs). If this is the case, then the analysis can focus on these quantities at the expense of others and reduce the number of forward calculations. In the literature this is termed as reducing the dimension of the problem, where “dimension” is related to the number of input/output quantities of interest (QoI) \cite{vohra2020,hou2022,xian2023,kim2024}.

In this paper, a surrogate model combined with a dimensionality reduction technique (DR) is established as a powerful inverse calculation method for uncertainty quantification (UQ). In this context, UQ can broadly be defined as the science of identifying, quantifying, and reducing uncertainties associated with predictive models, their outcomes, or quantities of interests. The approach is then applied to the complex inverse problem associated with a laterally loaded pile as an example, which has a high number of outputs of interest (i.e., high dimensionality). The method is shown to offer several compelling advantages in a sequential Bayesian inversion framework, namely:
\begin{itemize}[left=0pt]
 \item Scalability and adaptability: The DR-based surrogate is both scalable and adaptable, enabling it to address a wide range of engineering problems with high-dimensional outputs which are temporal in nature;
 \item Computational efficiency: A noteworthy feature is the significant reduction in computational complexity with respect to the output size $N$ and the number of required forward calculations $K$. This leads to a decrease in required computational resources;
 \item Improved prediction accuracy: With the ability to capture the dependency of output quantities, enhancements in the accuracy of model predictions are achieved;
\end{itemize}

This paper is structured as follows: \Cref{sec:ingredients} introduces the selected components for the DR-based surrogate model. \Cref{sec:DRSM} introduces the general formulation of the proposed DR-based surrogate and gives the computational details. \Cref{sec:numerics} tests the performance of the proposed method on a real pile problem and shows the superiority of the proposed approach. \Cref{sec:conclusion} gives brief conclusions. 

\section{Formulation of the DR-based surrogate modelling technique}\label{sec:ingredients}\noindent
Many of the most popular surrogate approaches, including Gaussian processes (kriging), polynomial chaos expansion (PCE), support vector machine (SVM), and neural networks (NN) \cite{torre2019,hou2022,vohra2020,kim2024,xian2023}, encounter challenges in parameterising and training surrogate models in high dimensional output space (large number QoIs) \cite{lataniotis2019} for the purpose of UQ. Similar challenges also arise in the presence of high-dimensional input space, often referred to as the \textit{curse of dimensionality} \cite{verleysen2005}. For high-dimensional outputs, current methods for constructing surrogate models are limited due to their reliance on constructing individual surrogate models, emulating individual QoI values \citep{lo2019,qi2017}. As many engineering problems involve high-dimensional, interrelated QoIs, this deficiency limits the applicability of many standard surrogate approaches. Problems in the civil engineering domain, for example, the responses of an offshore wind turbine under environmental loading, design of a retaining wall or construction of a tunnel in an urban environment, often lead to dependencies in the QoIs that may be overlooked in current surrogate model design.

Although methods such as multi-kriging \cite{toal2023} or multi-output regression \cite{rai2012} have the ability to overcome these shortcomings, the computational cost associated with using them to perform inverse analysis is demanding and often prohibitive. The reason for this is that as the output dimension increases (number of QoIs), the regression function needs to map the input variables to a larger number of QoIs. To efficiently manage high dimensions, an alternative approach is to employ DR techniques. Different flavours of DR techniques have been proposed to deal with high dimensional data, including linear compression such as \textit{principal component analysis} (PCA), nonlinear multiscaling (\textit{MDS}), kernel tricks exemplified by \textit{Kernal-PCA}, low dimensional embedding such as \textit{Isomap} and neural networks like \textit{autoencoder} (see \cite{kontolati2022,van2009} for a review).

DR-based surrogate modelling involves the fusion of two computational methodologies, which may include any methodology from the list of DR and surrogate modelling techniques referred to above. In this paper, the approach will be exemplified using PCA and PCE as the DR and surrogate modelling technique, respectively. This combination is used for simplicity and generality. However, other choices could still be used in the methodology described below. 
\subsection{DR methodology}\label{subsec: PCA}\noindent
Consider a computational model $\mathcal{M}$ (e.g., numerical analysis) with $K$ independent realisations of input space, defined as the \textit{training size}, $\mathcal{X} = \{\boldsymbol{x}^{(i)}\}_{i=1}^K \subset \mathcal{D}_{\boldsymbol{X}}$ such that certain QoIs, represented in a vector $\textbf{y} \in \mathbb{R}^N$, are a function of input parameters $\boldsymbol{x}$:
\begin{equation}
 \label{eq:UQ_model}
 \mathcal{M}:\boldsymbol{x} \in \mathcal{D}_{\boldsymbol{X}} \subseteq \mathbb{R}^M \mapsto
 \boldsymbol{\mathsf{y}} = \mathcal{M}(\boldsymbol{x}) \in \mathbb{R}^N 
\end{equation}

\noindent The set of computational outputs corresponding to the \textit{training size} can be collated into $\boldsymbol{Y} = \{\boldsymbol{\mathsf{y}}^{(i)}\}_{i=1}^K \subset \mathcal{D}_{\mathbf{Y}}$. This set can now be reduced following an appropriate DR technique. In an abstract form, the transformation from the \textit{original space} $\mathcal{D}_{\boldsymbol{Y}} \subseteq \mathbb{R}^{N}$ to a \textit{reduced space} $\mathcal{D}_{\boldsymbol{Z}} \subseteq \mathbb{R}^{N'}$ ($N' \ll N$) can be expressed as a DR mapping:
\begin{equation}
 \mathcal{T}_{DR} : \mathcal{D}_{\boldsymbol{Y}} \mapsto \mathcal{D}_{\boldsymbol{Z}}
\end{equation}
where the underlying assumption is that $\mathcal{D}_{\boldsymbol{Z}}$ is embedded inside $\mathcal{D}_{\boldsymbol{Y}}$. The nature and number of reduced outputs, $N'$, is dependent on the specific DR technique chosen. As noted above, this study adopted the simple but effective linear dimensionality reduction technique known as PCA \cite{vohra2020,nagel2020,wagner2020,wagner2021}. It is worth noting that PCA is also referred to as \textit{Karhunen-Loève expansion}, \textit{Hotelling
transform}, or \textit{proper orthogonal decomposition} depending on the specific scientific community.

PCA begins by calculating the estimation of the expectation $\boldsymbol{\mu_{Y}}$ and the covariance matrix $\boldsymbol{\Sigma_{Y}}$ of an independent and identically distributed (i.i.d) dataset:
\begin{equation}
\label{eq: PCA-mu}
\boldsymbol{\mu_{Y}} = \mathbb{E}[\boldsymbol{Y}]
= \begin{bmatrix}
\mu_{\boldsymbol{\mathsf{y}}_{1}} \\
\mu_{\boldsymbol{\mathsf{y}}_{2}} \\
\vdots \\
\mu_{\boldsymbol{\mathsf{y}}_{N}}
\end{bmatrix}^{\mathsf{T} }
\ 
\text{with}
\
\mu_{\boldsymbol{\mathsf{y}}_{i}} = \frac{1}{K} \sum_{k=1}^{K} \mathsf{y}_{i}^{k}
,\
i=1,\cdots,N
\end{equation}
\begin{equation} 
\label{eq: PCA-sigma}
\boldsymbol{\Sigma_{Y}} \approx \text{Cov}[\boldsymbol{Y}] 
 = \mathbb{E}[
(\boldsymbol{Y} -\boldsymbol{\mu_{Y}} )
(\boldsymbol{Y} - \boldsymbol{\mu_{Y}} )^{\mathsf{T}}
] 
\end{equation}
In which, the superscript denotes the forward realisation and the subscript denotes the output location. Since $\boldsymbol{\Sigma_{Y}}$ is symmetric and positive definite, one can find linearly independent eigenvectors $\boldsymbol{\phi}_{i}$ with positive eigenvalues $\lambda_{i}$ for $i=1,\cdots,N$. These characteristic vectors and values should satisfy:
\begin{equation} 
\label{eq: PCA-decompostion}
\boldsymbol{\Sigma_{Y}}\boldsymbol{\phi}_{i}
=
\lambda_{i}\boldsymbol{\phi}_{i}
\end{equation}
The $N$ eigenvectors of $\boldsymbol{\Sigma_{Y}}$ (new coordinates to be projected onto) are collected into column vectors $\boldsymbol{\Phi}_{N} = \{\boldsymbol{\phi}_{1},\cdots,\boldsymbol{\phi}_{N}\}$ for $i=1,\cdots,N$. The corresponding eigenvalue $\lambda_{i}$ signifies the variance of $\boldsymbol{Y}$ in the direction of the $i \text{-th}$ principle component showing the descending order $\lambda_{1}\ge \lambda_{2} \ge \cdots \ge \lambda_{N}$. Since the original data $\boldsymbol{Y}$ has now been centred and decorrelated in \Cref{eq: PCA-mu,eq: PCA-sigma,eq: PCA-decompostion}, the linearly transformed random vectors have a zero mean and a diagonal covariance matrix. Orthogonal components can be then fully represented as:
\begin{equation} 
\boldsymbol{Z}_{\text{full}} = \boldsymbol{\Phi}_{N}^{\mathsf{T}}
(\boldsymbol{Y} - \boldsymbol{\mu_{Y}})
\end{equation}
If the QoIs have some correlation with each other, then only a subset of the principal components of $\boldsymbol{Z}_{\text{full}}$ are needed to obtain an accurate approximation. As the influence of each principal component reduces as $i$ increases from $1$ to $N$, then this implies that only a reduced number of principal values, $N'$, are required.

\begin{figure}[htbp]
 \centering 
 \begin{subfigure}[b]{0.8\textwidth}
 \centering
 \includegraphics[width=50mm]{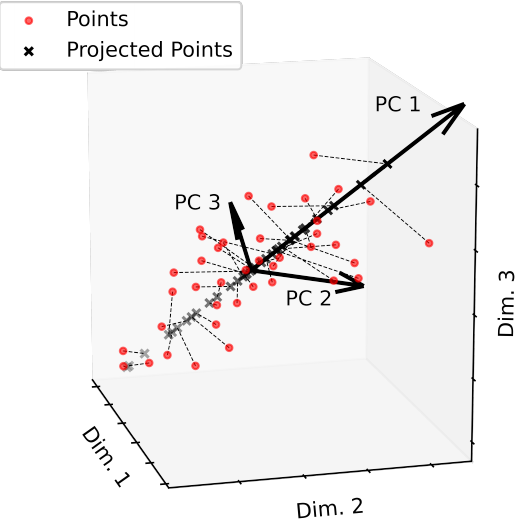}
 \caption{\textit{Decomposition of a 3D data to PC1}}
 \label{fig: PCA3D}
 \end{subfigure}
 \begin{subfigure}[b]{0.8\textwidth}
 \centering
 \includegraphics[width=90mm]{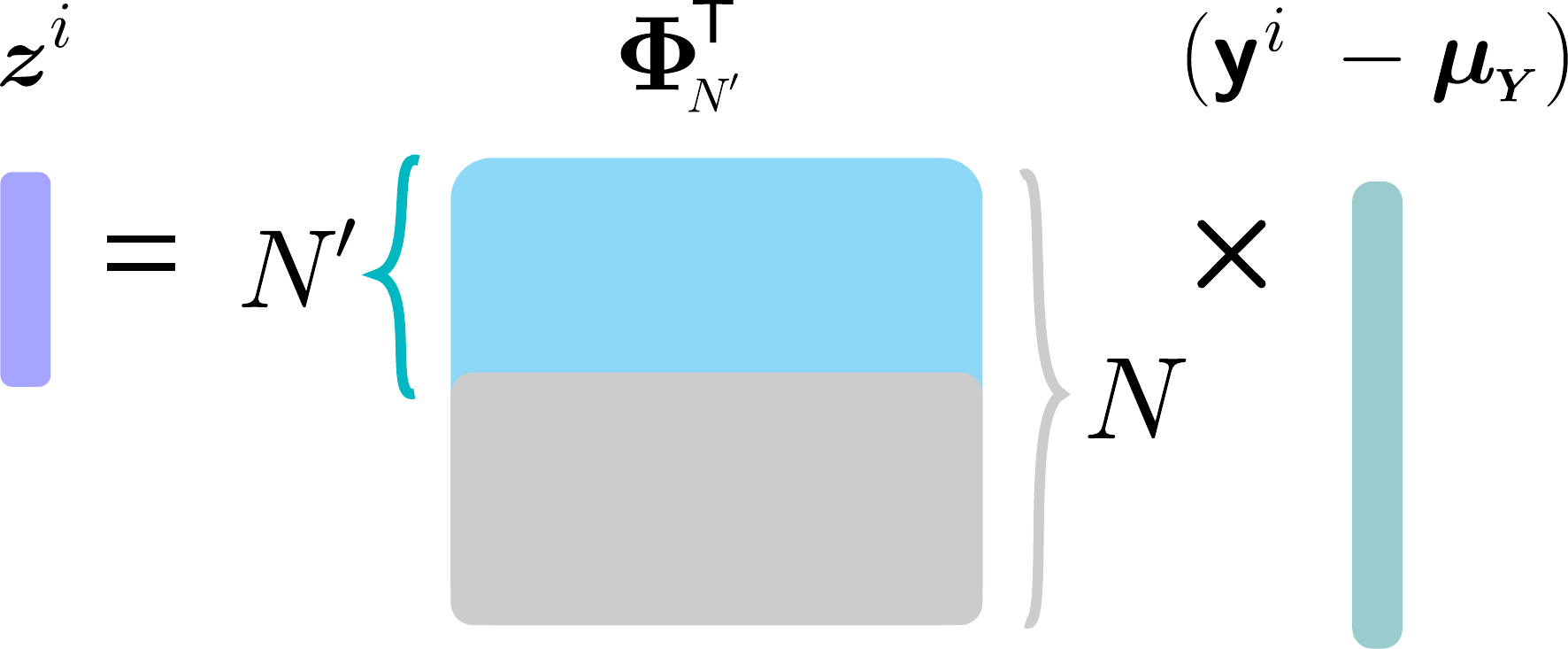}
 \caption{\textit{One realisation of PCA with selected eigenvectors for N-dimensional data}}
 \label{fig: PCA-schematic}
 \end{subfigure}
 \caption{\textit{A visualised PCA process}}
\end{figure}

 By retaining only $N'$ principal components (PCs) with the highest variance, the model outputs $\boldsymbol{Y}$ can be approximated as $\boldsymbol{Z}$:
 \begin{equation}
 \label{eq: PCA-component}
\boldsymbol{Z} = \boldsymbol{\Phi}^{\mathsf{T}}_{N'}
(\boldsymbol{Y} - \boldsymbol{\mu_{Y}}) \approx \boldsymbol{Z}_{\text{full}}
\end{equation} 
Principles of decomposition are shown in \Cref{fig: PCA3D} and the output vectors of each realisation (\textit{observation sample}) can be reduced through \Cref{eq: PCA-component}. A specific realisation $\boldsymbol{x}^{i}$ is visualised in \Cref{fig: PCA-schematic} based on the selected eigenvectors $ \boldsymbol{\Phi}^{\mathsf{T}}_{N'}$. 
Subsequently, the outputs of all realisations can be compressed while retaining most of the total variations by:
\begin{equation}
 \label{PCA-components_expression_Z}
 \boldsymbol{z} = 
\begin{bmatrix}
({\boldsymbol{z}^{1}})^{\mathsf{T}} \\
({\boldsymbol{z}^{2}})^{\mathsf{T}} \\
\vdots \\
({\boldsymbol{z}^{K}})^{\mathsf{T}} 
\end{bmatrix}
=
\begin{bmatrix}
z_{1}^{1} & z_{2}^{1} & \cdots & z_{N'}^{1} \\
z_{1}^{2} & z_{2}^{2} & \cdots & z_{N'}^{2} \\
\vdots & \vdots & \ddots & \vdots \\
z_{1}^{K} & z_{2}^{K} & \cdots & z_{N'}^{K}
\end{bmatrix} 
\end{equation}

The selection of the number $N'$ can be chosen by specifying $\sum_{i=1}^{N'} \lambda_{i} = 
(1-\varepsilon_{\text{DR}}^{threshold})\sum_{i=1}^{N} \lambda_{i}$, where $\varepsilon_{\text{DR}}^{threshold}$ is a predefined threshold. An adequate number of principal components should be selected to represent the system in an optimal manner. If too few PCs are selected, a poor model will be obtained. Conversely, if too many PCs are selected, negligible improvement in computation costs will be gained. In high-dimensional data outputs (many QoIs), over-parameterisation may also inadvertently capture or introduce unwanted noise in the reconstructed dataset. To avoid these problems, several criteria for selecting the optimum number of PCs have been proposed, such as scree plot, explain-variance, permutation test, cross-validation and variance of reconstruction error \cite{valle1999,saccenti2015,qin2000,mnassri2010}. With a suitable threshold and selection criteria, the original space can be reconstructed as $\boldsymbol{Y}^{Re}$ through its optimum $N'$ principal components using
\begin{equation}
\label{eq: PCA-reconstruct}
\boldsymbol{Y} 
=\boldsymbol{\mu_{Y}} + 
\sum_{i=1}^{N} \boldsymbol{z}_{i}\boldsymbol{\phi}_{i}
\approx \boldsymbol{Y}^{Re} 
= \boldsymbol{\mu_{Y}} + 
\sum_{i=1}^{N'} \boldsymbol{z}_{i}\boldsymbol{\phi}_{i}
\end{equation}

\begin{figure}[htbp]
 \centering \includegraphics[width=90mm]{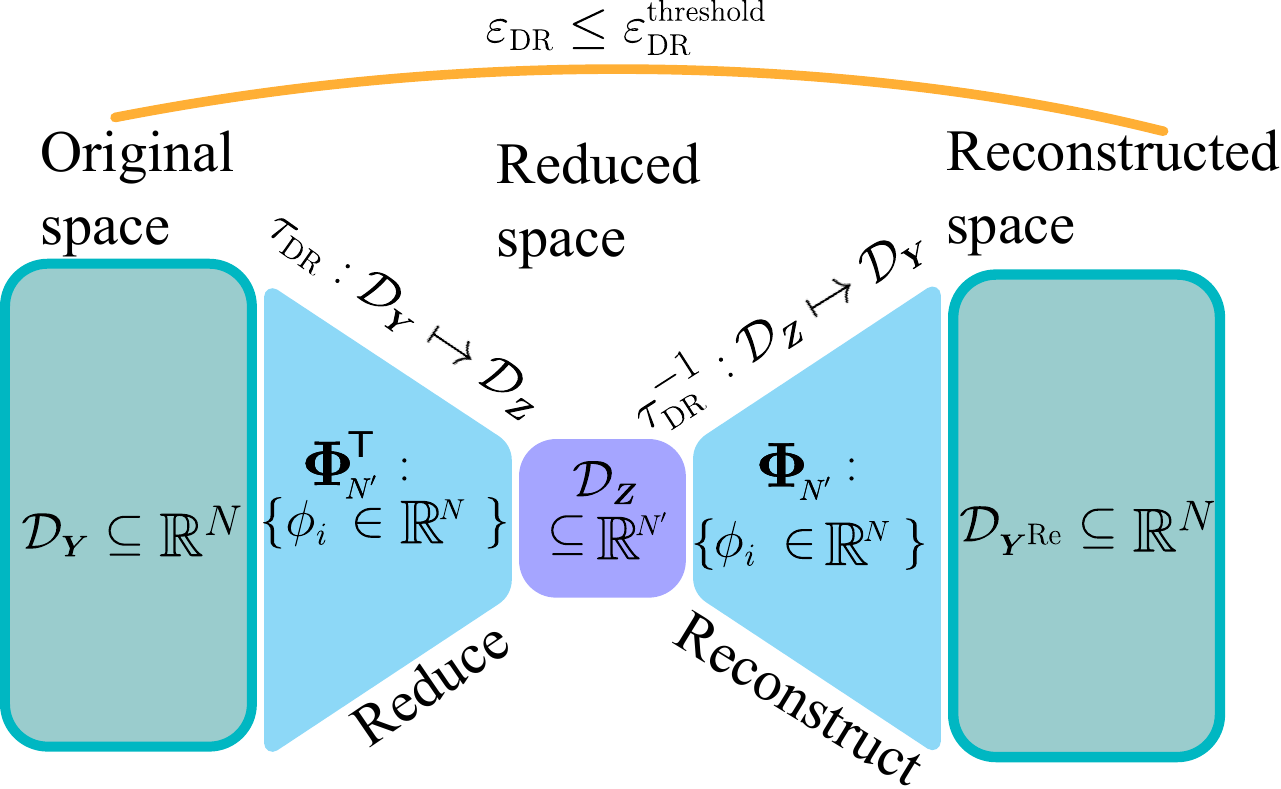}
 \caption{\textit{DR-flowchart}}
 \label{fig: PCA-flowchart}
\end{figure}

Generally, most DR techniques aim to find a lower dimensional representation of the dataset while minimising the loss of accuracy due to compression, as illustrated in \Cref{fig: PCA-flowchart}. All methods seek to (1) calculate the reduced spaces $\boldsymbol{Z} $ (e.g., principal components in PCA) and (2) based on a predefined threshold $\varepsilon_{\text{DR}}^{threshold}$ and transformation process $\mathcal{T}_{DR}^{-1}$, reconstruct the outputs $\boldsymbol{Y}^{Re}$. While PCA is a highly effective and powerful compression tool, the main idea of PCA is still based on linear decomposition. As data complexity increases, the limitations of PCA become apparent. In such cases, more advanced DR techniques, such as \textit{MDS}, \textit{kPCA}, and \textit{autoencoder}, should be considered.

\subsection{Surrogate model}\label{subsec: PCE}\noindent
Polynomial chaos expansions (PCE) is a commonly used spectral method for approximating functions in place of an expensive model and has been used extensively in engineering practice \cite{lataniotis2019,Breitenmoser2023,wagner2020} for UQ. Among the many available classes of surrogate models, PCE stands out for its simple polynomial form and ease of deployment \cite{torre2019}. As a non-intrusive surrogate method, PCE aims at providing a functional approximation of a forward model through its spectral representation on a suitably built basis of polynomial functions. To build a PCE, consider a random vector input $\boldsymbol{X} \in \mathbb{R}^{M}$ characterised by the joint probability density function (PDF) $f_{\boldsymbol{X}}$, and a computation model with finite variance satisfying:
\begin{equation}
 \mathbb{E}[\boldsymbol{Y}^2] 
 =\int_{\mathcal{D}_{\boldsymbol{X}}}
 \mathcal{M}^2(\boldsymbol{x})f_{\boldsymbol{X}}(\boldsymbol{x})
 d\boldsymbol{x}< \infty 
\end{equation}
a \textit{PCE} of the forward model can then be represented as :
\begin{equation}
\label{eq: PCE_basis}
\textbf{y}=\mathcal{M}(\boldsymbol{x}) =
\sum\limits_{\alpha \in \mathbb{N}^M }{\boldsymbol{c}_{\boldsymbol{\alpha}} \Psi_{\boldsymbol{\alpha}}(\boldsymbol{x})}
\end{equation}
where $\Psi_{\boldsymbol{\alpha}}(\boldsymbol{x})$ are multivariate polynomials orthonormal relative to $f_{\boldsymbol{X}}$, $\boldsymbol{\alpha} \in \mathbb{N}^M $ identifies the multi-indices storing the degrees of the multivariate polynomials $\Psi_{\boldsymbol{\alpha}}$, and $\boldsymbol{c}_{\boldsymbol{\alpha}} \in \mathbb{R}$ are the corresponding coefficients.

The polynomial basis $\Psi_{\boldsymbol{\alpha}}(\boldsymbol{x})$ in \Cref{eq: PCE_basis} is conventionally constructed using a set of \textit{univariate orthonormal polynomials} $\psi_{k}^{i}(\boldsymbol{x}_{i})$, where these polynomials satisfy:
\begin{equation}
 \left \langle 
\psi_{j}^{i},\psi_{k}^{i}
 \right \rangle 
= \delta_{jk}
\end{equation}
where $i$ identifies the input variable and $\delta_{jk}$ is the \textit{Kronecker symbol}. By taking the tensor product of their univariate counterparts, $\Psi_{\alpha}$ can be represented as:
\begin{equation}
\Psi_{\boldsymbol{\alpha}}(\boldsymbol{x})
 \overset{\mathrm{def}}{=}
\prod_{i=1}^{M} 
\psi_{\alpha_{i}}^{i}(x_{i})
\end{equation}
in which, $\boldsymbol{\alpha}= ( \alpha_{1},\cdots,\alpha_{M} ) \in \mathbb{N}^{M}$, and $\Psi_{\boldsymbol{\alpha}}(\boldsymbol{x})$ satisfies orthonormal properties as:
$ \left \langle 
\Psi_{\boldsymbol{\alpha}}(\boldsymbol{x}),\Psi_{\boldsymbol{\beta}}(\boldsymbol{x})
 \right \rangle 
= \delta_{\boldsymbol{\alpha}\boldsymbol{\beta}}$, where the symbol $\delta_{\boldsymbol{\alpha}\boldsymbol{\beta}}$ is an extension of the \textit{Kronecker symbol} to the multi-dimensional case.

The infinite series expansion in \Cref{eq: PCE_basis} cannot be handled in realistic applications. Consequently, a truncation scheme is introduced, e.g., by selecting the maximal degree of the polynomials. The coefficients vector, $\boldsymbol{c}_{\boldsymbol{\alpha}}$, in \Cref{eq: PCE_basis} is typically estimated by a least-squares analysis. However, in this study, due to the computational costs associated with the forward model and high-dimensional nature of the output field, the sparse PCE Least Angle Regression algorithm has been used (see \cite{blatman2011} for more detail). After calculating the coefficients of the PCE for a given basis, the truncated \textit{mulivariate surrogate PCE}, can then be written as
\begin{equation}
\label{eq: surrogatePCE}
\mathcal{M}(\boldsymbol{x}) \approx \mathcal{M}^{PC}(\boldsymbol{x})=
\sum\limits_{\boldsymbol{\alpha} \in \mathcal{A} }{\boldsymbol{c}_{\boldsymbol{\alpha}} \Psi_{\boldsymbol{\alpha}}(\boldsymbol{x})}
\end{equation}
\noindent where $\mathcal{A}\subset\mathbb{N}^{M}$ is the set of selected multi-indices of the univariate polynomials. 

To assess the accuracy of the obtained PCE, error estimators are essential for quantifying the fidelity of a surrogate model in approximating the original model. To determine the precision of a surrogate model, \textit{cross validation} or \textit{leave-one-out test} can be used to assess the error. Once the coefficients vector $\boldsymbol{c}_{\boldsymbol{\alpha}}$ and the surrogate are constructed, explicit expressions for moments characterising the model outputs can be extracted. $\tilde{\mathcal{M}}(\boldsymbol{x})$ can be subsequently calculated at a negligible cost compared to the original forward model.

\subsection{The proposed DR-based surrogate: PCA-PCE}
\begin{figure}[!ht]
 \centering
 \includegraphics[width=90mm]{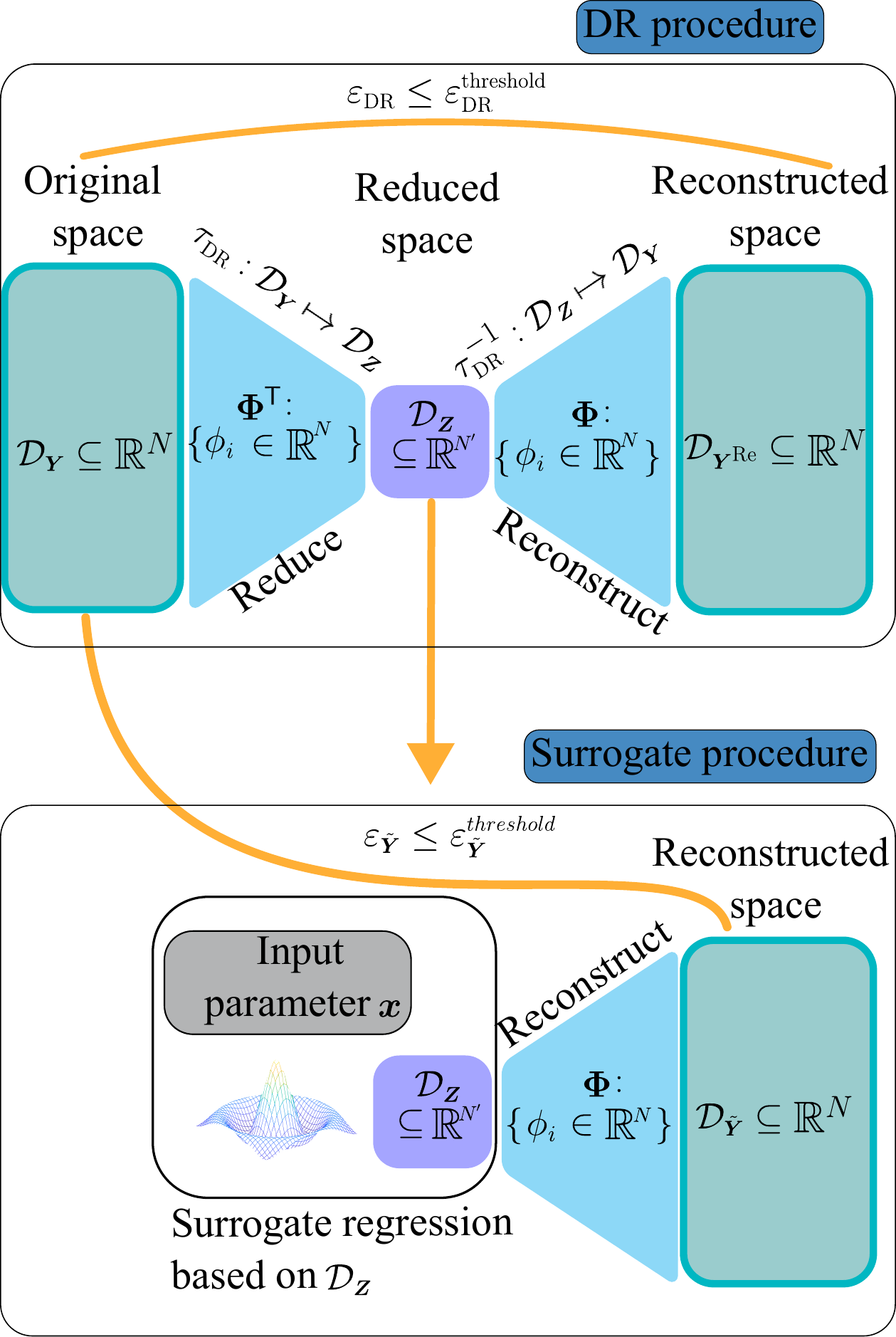}
 \caption{A two-step surrogate model construction: PCA-PCE}
 \label{fig: PCA-PCE}
\end{figure}
\noindent Combining PCA (\Cref{subsec: PCA}) and PCE (\Cref{subsec: PCE})
constitutes a particularly efficient surrogate modelling technique. Starting from realisations of input space $\boldsymbol{X}$ and corresponding high-dimensional model outputs, $\boldsymbol{Y}$, the principal component transformation for the output space can be readily derived with \Cref{eq: PCA-component}. Each component $\boldsymbol{z}_{i}$ is then expanded onto a polynomial basis using \Cref{eq: PCE_basis} such that
\begin{equation}
\boldsymbol{z}_{i} \approx
\sum\limits_{\boldsymbol{\alpha} \in \mathcal{A} }{\boldsymbol{c}_{i,\boldsymbol{\alpha}} \Psi_{\boldsymbol{\alpha}}(\boldsymbol{X})}
\end{equation}
This expression together with \Cref{eq: PCA-reconstruct} yields a surrogate model called PCA-PCE that relates the input to the high-dimensional model outputs by:
\begin{equation}
\label{eq: PCA-PCE}
\boldsymbol{Y} 
\approx \boldsymbol{\mu_{Y}} + 
\sum_{i=1}^{N'} \boldsymbol{z}_{i}\boldsymbol{\phi}_{i}
\approx \boldsymbol{\mu_{Y}} + 
\sum_{i=1}^{N'} (\sum\limits_{\boldsymbol{\alpha} \in \mathcal{A} }{\boldsymbol{c}_{i,\boldsymbol{\alpha}} \Psi_{\boldsymbol{\alpha}}(\boldsymbol{X})}) \boldsymbol{\phi}_{i}
\end{equation}

The combination of both approaches can together be used as a new class of surrogate model which is able to approximate high-dimensional outputs with greater efficiency and accuracy. Based on this, a two-step DR-surrogate can be constructed as outlined in \Cref{fig: PCA-PCE} and as summarised below:
\begin{itemize}[left = 0pt]
 \item DR procedure: a PCA step is carried out to extract data features $\boldsymbol{\Phi}^{\mathsf{T}}$ of the original high-dimensional output space and compress the outputs $\mathcal{D}_{\boldsymbol{Y}}$ into the reduced space $\mathcal{D}_{\boldsymbol{Z}}$. Based on a predefined reconstruction error $\varepsilon_{\text{DR}}^{threshold}$, the optimal reduced principal number $N'$ is selected.
 \item Surrogate procedure: Based on the obtained reduced space from the DR procedure, a surrogate (e.g., PCE) $\tilde{\mathcal{M}}$ is constructed directly. To test the surrogate built on the reduced space, model responses are reconstructed into the original space based on the prescribed $\varepsilon_{\tilde{\boldsymbol{Y}}}^{threshold}$.
\end{itemize}
in which, $\mathcal{D}_{\boldsymbol{Y}^{Re}}$ denotes the data space for PCA reconstruction, and $\mathcal{D}_{\tilde{\boldsymbol{Y}}}$ denotes the data space for PCA-PCE reconstruction.

\section{Inverse analysis}\label{sec:DRSM}\noindent
\subsection{Bayesian inversion}\noindent
All computational models (i.e., numerical analyses) are simplifications of reality, and field/experimental observations will contain some level of \textit{uncertainty} with the computer model. To explicitly link this mismatch, one option is to model the \textit{uncertainty} as an \textit{additive} observation error $\boldsymbol{\varepsilon}$ between the model prediction $\mathcal{M}(\boldsymbol{x})$ and the observations $\boldsymbol{y}$:
\begin{equation}
 \label{eq: modelling_discrepancy}
 \boldsymbol{y} = \mathcal{M}(\boldsymbol{x}) + \boldsymbol{\varepsilon}
\end{equation}
\noindent where $\boldsymbol{\varepsilon}\in\mathbb{R}^N$ is the uncertainty that encapsulates all the discrepancies between the field/experimental observation $\boldsymbol{y}$ and the model prediction. The collection of $m$ sets of observations can be gathered into a set given by:
\begin{equation}
\label{eq:high-outputs}
 \mathcal{Y} = \{\boldsymbol{y}^{(i)}\}_{i=1}^m = 
\{[ y_{1}^{i}, y_{2}^{i} , \cdots , y_{N}^{i} ]\}_{i=1}^m = \{ \mathcal{M}(\boldsymbol{x}) +\boldsymbol{\varepsilon}^{(i)}\}_{i=1}^m
\end{equation}

Within the Bayesian framework, a general solution to an inverse problem is the posterior $\pi(\boldsymbol{x}|\mathcal{Y})$ defined over the input domain $\mathcal{D}_{\boldsymbol{X}}$, formally expressed using \textit{Bayes's theorem}:
\begin{equation}
\label{eq: full_Bayes_theorem}
\pi(\boldsymbol{x}|\mathcal{Y}) = \frac{{\mathcal{L}(\boldsymbol{x}|\mathcal{Y}) \cdot \pi(\boldsymbol{x})}}{{\pi(\mathcal{Y})}} 
\end{equation}
where $\pi(\boldsymbol{x}|\mathcal{Y})$ is calculated from the prior $\pi(\boldsymbol{x})$, likelihood $\mathcal{L}(\boldsymbol{x}|\mathcal{Y})\stackrel{\mathrm{def}}{=}\pi(\mathcal{Y}|\boldsymbol{x})$
and the evidence $\pi(\mathcal{Y})$. 

\Cref{eq: full_Bayes_theorem} is rarely used directly, as computing the evidence $\pi(\mathcal{Y})$ is typically intractable due to model complexity or high computational cost. Thus, methods which can approximate the posterior distribution are used. One common approach is to directly sample from the posterior distribution using a Markov Chain Monte Carlo (MCMC) sampling algorithm. To avoid the computational burden of sampling using the forward model, MCMC simulations use the PCA-PCE surrogate to construct the Markov chains ($\boldsymbol{X}^{(1)}, \boldsymbol{X}^{(2)}, \cdots )$ over the prior support $\mathcal{D}_{\boldsymbol{X}}$. In this study, the \textit{affine-invariant ensemble sampler} (AIES) algorithm \cite{wagner2020} was used to construct the Markov chains with a burn-in period of 70\%.

To verify the calibration process and inversion results, predictive distributions of the model need to be calculated. One predictive distribution of interest is the \textit{prior predictive distribution} $\pi(\boldsymbol{\mathsf{y}})$ which is defined as:
\begin{equation}
 \label{eq: prior predictive}
 \pi(\boldsymbol{\mathsf{y}}) = \int_{\mathcal{D}_{\boldsymbol{X}}} 
 \pi(\boldsymbol{x}) \pi(\boldsymbol{\mathsf{y}}|\boldsymbol{x}) {\rm{d}} \boldsymbol{x}
\end{equation}
This distribution expresses beliefs about the future and summarises the uncertainties before the calibration. It is essential for identifying and ruling out challenging inverse problems before proceeding with expensive calibrations. In practical terms, it helps determine if the field/ experimental observations can be reproduced by assessing whether they fall within the predictive distribution. The \textit{posterior predictive} distribution $\pi(\boldsymbol{\mathsf{y}}|\mathcal{Y})$ summarises the uncertainty of the input vector after the calibration and can be expressed in the form:
\begin{equation}
 \label{eq: posterior predictive}
 \pi(\boldsymbol{\mathsf{y}}|\mathcal{Y}) = \int_{\mathcal{D}_{\boldsymbol{X}|\mathcal{Y}}} 
 \pi(\boldsymbol{x}|\mathcal{Y}) \pi(\boldsymbol{\mathsf{y}}|\boldsymbol{x}) {\rm{d}} \boldsymbol{x}
\end{equation}
Within the $\mathcal{D}_{\boldsymbol{X}|\mathcal{Y}}$ domain, the parameter set which is located at the maximum value of the posterior (\textit{maximum a posterior}, MAP ) is also shown along with the predictive distributions by:
\begin{equation}
 \label{eq: MAP}
 \begin{aligned}
 \boldsymbol{x}^{\rm{MAP}} &= \mathop{\arg\max}\limits_{\boldsymbol{x} \in \mathcal{D}_{\boldsymbol{X}}}
 \pi(\boldsymbol{x}|\mathcal{Y}) \\
 &=\mathop{\arg\max}\limits_{\boldsymbol{x} \in \mathcal{D}_{\boldsymbol{X}}}
 {\mathcal{L}(\boldsymbol{x}|\mathcal{Y}) \pi(\boldsymbol{x})} 
 \end{aligned}
\end{equation}

\subsection{Incorporating uncertainties into sequential Bayesian inversion}\label{subsec: modelling_error}\noindent
This study assumes that the \textit{uncertainty} term $\boldsymbol{\varepsilon}$ follows a zero mean multivariate Gaussian distribution, i.e., $\boldsymbol{\varepsilon} \in \mathcal{N}(\boldsymbol{\varepsilon}|\boldsymbol{0},\boldsymbol{\Sigma})$. Given $m$ sets of sparse independent observation vectors $\mathcal{Y} {=} (\boldsymbol{y}_1,\cdots,\boldsymbol{y}_m)$ and following the construction of an appropriate surrogate $\tilde{\mathcal{M}}(\boldsymbol{x})$, the likelihood can be explicitly expressed as:
\begin{equation} 
 \label{eq: modified_Likelihood function}
\begin{aligned}
 \mathcal{L}(\boldsymbol{x}|\mathcal{Y}) =& \prod_{i=1}^{m} N(\boldsymbol{y}_i|\tilde{\mathcal{M}}(\boldsymbol{x}),\boldsymbol{\Sigma}) \\
 =& \prod_{i=1}^{m}\frac{1}{\sqrt{(2 \pi)^{N}{\rm{det}} 
 (\boldsymbol{\Sigma})}}\exp\left(-\frac{1}{2}\left[\boldsymbol{y}_i - \tilde{\mathcal{M}}(\boldsymbol{x})\right]^{\mathsf{T}} \boldsymbol{\Sigma}^{-1}\left[\boldsymbol{y}_i - \tilde{\mathcal{M}}(\boldsymbol{x})\right]\right)
\end{aligned}
\end{equation} 

The exact residual covariance matrix, $\boldsymbol{\Sigma}$, is most commonly unknown. However, by parameterising the matrix as $\boldsymbol{\Sigma}(\boldsymbol{x^{\varepsilon}})$, one may treat its parameter $\boldsymbol{x^{\varepsilon}}$ as additional unknowns that can be inferred jointly with the input parameters of $\mathcal{M}$. In this setting, the parameter vector is defined by $\boldsymbol{x} = (\boldsymbol{x}^{\mathcal{M}},\boldsymbol{x^{\varepsilon}})$, i.e., a combined vector of \textit{forward model parameters} $\boldsymbol{x}^{\mathcal{M}}$ and \textit{uncertainty parameters} $\boldsymbol{x^{\varepsilon}}$. To characterise $\boldsymbol{x^{\varepsilon}}$, a diagonal covariance matrix of the form $\boldsymbol{\Sigma} = \sigma^2 \times \boldsymbol{I}$ with an unknown residual variance $\sigma^2$ is adopted, where $\boldsymbol{I}$ denotes an identity matrix. In this way, the \textit{uncertainty parameters} $\boldsymbol{x^{\varepsilon}}$ can be reduced to a single scalar, i.e., $\boldsymbol{x}^{\varepsilon} \equiv \sigma^2$. By treating $\boldsymbol{x}^{\mathcal{M}}$ and $\boldsymbol{x^{\varepsilon}}$ as priorly independent, the joint prior distribution can be represented as:
\begin{equation}
 \label{eq: joint_prior}
\pi(\boldsymbol{x}) = \pi(\boldsymbol{x}^{\mathcal{M}})\pi(\sigma^2)
\end{equation}
and the corresponding posterior distribution can be computed as:
\begin{equation}
 \label{eq: posterior_final}
\begin{aligned}
\pi(\boldsymbol{x}^{\mathcal{M}},\sigma^2|\mathcal{Y} )
 &= \frac{\pi(\boldsymbol{x}^{\mathcal{M}})\pi(\sigma^2)
\mathcal{L}(\boldsymbol{x}^{\mathcal{M}},\sigma^2;\mathcal{Y} )}{\pi(\mathcal{Y})} \\
&\propto\pi(\boldsymbol{x}^{\mathcal{M}})\pi(\sigma^2)
\mathcal{L}(\boldsymbol{x}^{\mathcal{M}},\sigma^2;\mathcal{Y} )
\end{aligned}
\end{equation}

\begin{figure}[htbp]
\centering
 \includegraphics[width=140mm]{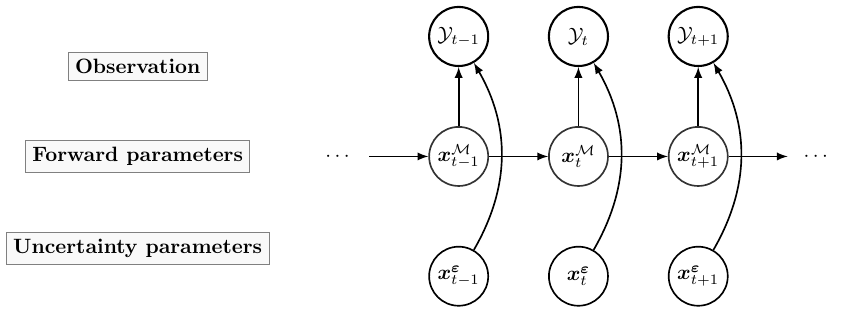}
 \caption{Sequential Bayesian inversion}
 \label{fig: sequential Bayesian inference}
\end{figure}

A virtue of the Bayesian framework is its ability to update parameter vectors $(\boldsymbol{x}^{\mathcal{M}},\boldsymbol{x}^{\varepsilon})$ when new observations $\mathcal{Y}_t$ become available. This process, known as data assimilation, allows the model to incorporate the new information. This feature is very useful for staged engineering projects as it enables different information flows for $\boldsymbol{x}^{\mathcal{M}}_{t}$ and $\boldsymbol{x^{\varepsilon}}_{t}$ at given stages $t$. By assuming each stage is only dependent on the preceding stage, the problem can be re-interpreted as a Markovian process in which Bayesian inversion is being sequentially applied. As shown in \Cref{fig: sequential Bayesian inference}, when new observations $\mathcal{Y}_{t}$ are obtained, the current posterior expression (otherwise known as the belief state) can be updated as:
\begin{equation}
 \label{eq: SBI_conditional_probability}
\pi(\boldsymbol{x}^{\mathcal{M}}_{t},\boldsymbol{x^{\varepsilon}}_{t}|\mathcal{Y}_{1:t})
\propto
\pi(\mathcal{Y}_{t}|\boldsymbol{x}^{\mathcal{M}}_{t},\boldsymbol{x^{\varepsilon}}_{t}) 
\pi(\boldsymbol{x}_{t}|\mathcal{Y}_{1:t-1})
\end{equation}
In this sequential inversion framework, only the \textit{forward model parameters} $\boldsymbol{x}^{\mathcal{M}}_{t}$ were updated and propagated to subsequent inversion stages. While the \textit{uncertainty parameters} $\boldsymbol{x^{\varepsilon}}_{t}$ were included in the inversion alongside $\boldsymbol{x}^{\mathcal{M}}_{t}$, they are not carried forward to subsequent stages. The rationale behind this lies in the fact that the framework considers each stage as a discrete event whereby the probability of the current true state is conditionally independent of the other earlier states. For prediction/control purposes, the belief state has proven to be a sufficient statistic \cite{aastrom1965}, i.e., the Markovian process is valid. As such, the information passing and updating process can be expressed as \Cref{eq: SSM_predict} and \Cref{eq: SSM_correct}, respectively:
\begin{equation}
 \label{eq: SSM_predict}
\pi(\boldsymbol{x}^{\mathcal{M}}_{t}|\boldsymbol{x}^{\mathcal{M}}_{0:t-1})
\propto
\pi(\boldsymbol{x}^{\mathcal{M}}_{t}|\boldsymbol{x}^{\mathcal{M}}_{t-1})
\end{equation}
\begin{equation}
 \label{eq: SSM_correct}
\pi(\mathcal{Y}_{t}|\boldsymbol{x}^{\mathcal{M}}_{0:t},\boldsymbol{x^{\varepsilon}}_{1:t}) 
\propto
\pi(\mathcal{Y}_{t}|\boldsymbol{x}^{\mathcal{M}}_{t},\boldsymbol{x^{\varepsilon}}_{t}) 
\end{equation}
where the superscript denotes the parameter types and the subscript denotes the time stages. 

For the priors shown in \Cref{eq: joint_prior,eq: posterior_final,eq: SBI_conditional_probability,eq: SSM_predict,eq: SSM_correct}, $\pi(\boldsymbol{x}^{\mathcal{M}})$ can be constructed based on existing heuristic knowledge or laboratory data prior to the collection of real-world field/experimental observations, and the calculated posterior can be passed as a prior for the next stage. $\pi(\sigma^2)$ is assumed to be stage-dependent and related to the collected observations $\mathcal{Y}$. Specifically, in this study, $\sigma^2$ is assumed as the maximum of the observation $\mathcal{Y}$ following $\sigma^2 \sim \mathcal{U}(0,\text{max}(\mathcal{Y}))$.

Notably, the updating process is not limited to a specific type of uncertainty. Alongside observational error $\boldsymbol{\varepsilon}$, more complex uncertainties, such as truncation errors $\eta$ related to surrogates $\tilde{\mathcal{M}}$, as well as model discrepancies $\boldsymbol{\delta}$, can be also considered as shown in \Cref{eq: COV_uncertainties}. Assuming these error terms also follow Gaussian distributions, the errors can be seamlessly integrated into a total error covariance matrix, represented as $\boldsymbol{\Sigma} = (\boldsymbol{\Sigma}_{T} + \boldsymbol{\Sigma}_{G} + \boldsymbol{\Sigma}_{D})$. Here, $\boldsymbol{\Sigma}_{T}$ represents the covariance of the surrogate truncation error, $\boldsymbol{\Sigma}_{G}$ denotes the covariance of the model discrepancy error, and $\boldsymbol{\Sigma}_{D}$ signifies the covariance of the observation error. 

\begin{equation}
\label{eq: COV_uncertainties}
\begin{aligned}
 \mathcal{Y}
=
\zeta(\boldsymbol{x}) +\boldsymbol{\varepsilon}
= &&
{\mathcal{M}}(\boldsymbol{x}) + {\delta(\boldsymbol{x})}
+ \boldsymbol{\varepsilon} \\
= &&
\tilde{\mathcal{M}}(\boldsymbol{x})
+ \eta({\boldsymbol{x}})
+ \delta(\boldsymbol{x})
+ \boldsymbol{\varepsilon}
\end{aligned}
\end{equation}

\subsection{Incorporating PCA-PCE into Bayesian inversion}\noindent\label{sec: PCA-PCE steps}
The quality of the inversion analysis in computing the posterior distribution $\pi(\boldsymbol{x}|\mathcal{Y})$ is directly affected by the efficacy and accuracy of the constructed surrogate. As described earlier, to obtain reliable predictive estimates for problems which involve high-dimensional outputs, the surrogate model must be capable of generating plausible outputs. A key limitation of traditional surrogate models is their inability to \textit{explicitly} capture connections among outputs. This often contradicts reality, as many QoIs exhibit temporal or spatial correlation. The proposed PCA-PCE surrogate modelling approach offers a promising solution to this issue. It combines the power of PCA to simplify high-dimensional outputs into a separate, lower-dimensional space with the efficiency of PCEs in modelling multivariate outputs. This approach enables a new type of surrogate that conducts MCMC simulations on a specific subspace rather than the entire output space. This approach therefore reduces the complexity of the inversion problem while preserving the inherent covariance structure of the original outputs.

\begin{figure}[htbp]
\centering
 \includegraphics[width=140mm]{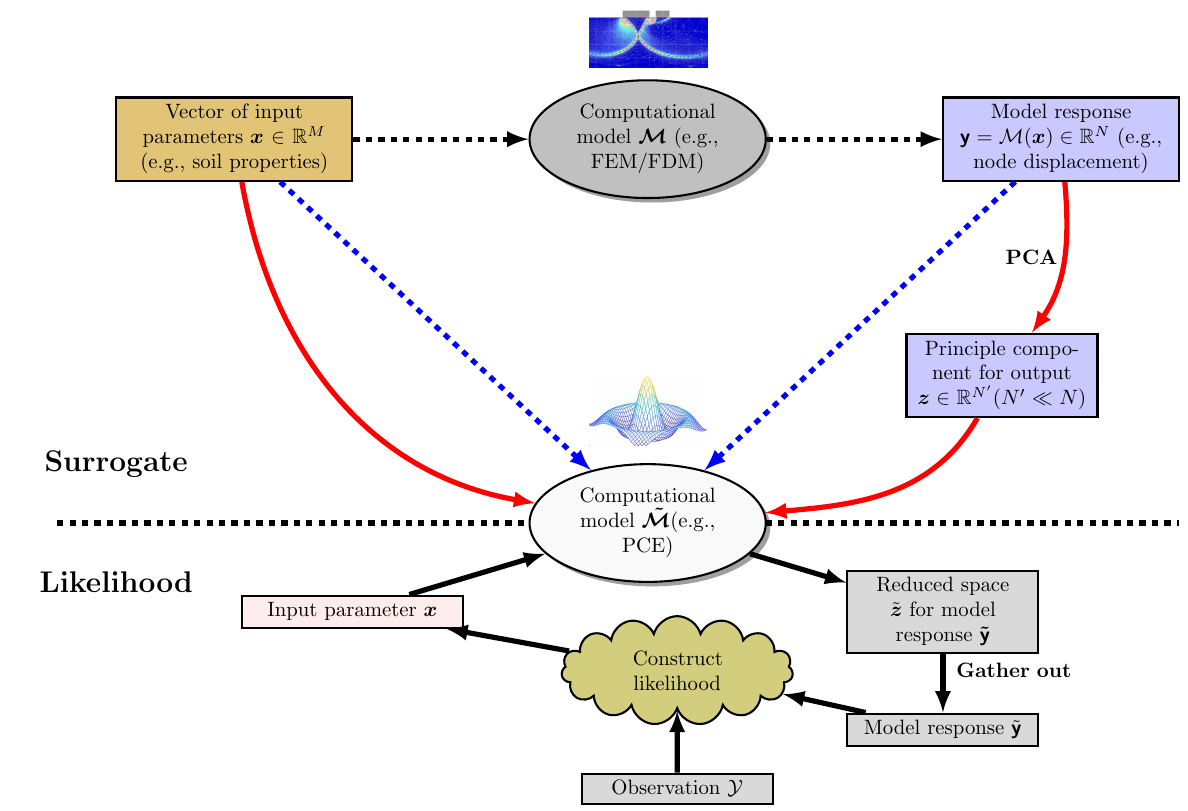}
 \caption{PCA-PCE surrogate in Bayesian inversion}
 \label{fig: PCA_PCE_likelihood}
\end{figure}

During the inversion process, the likelihood function $\mathcal{L}(\boldsymbol{x}|\mathcal{Y})$ is repeatably evaluated to adjust the uncertainties in the input space. Since observations $\mathcal{Y}$ hold physical meaning and the surrogate model is constructed within a non-physical reduced space, direct comparison between these two components is not meaningful. Therefore, an invertible mapping is necessary to transform reduced space to the original space for the likelihood computation. A brief step-by-step flowchart of the PCA-PCE surrogate in Bayesian framework is given in \Cref{fig: PCA_PCE_likelihood} as following: 
\begin{itemize}[left = 0cm]
 \item Step one: Based on the initial priors $\pi(\boldsymbol{x}^{\mathcal{M}})$, conduct FE modelling $\mathcal{X}$ using $K$ sets of samples. Evaluate the forward model $\mathcal{M}(\boldsymbol{x})$ at $\mathcal{X}$ and store the high-dimensional outputs $\boldsymbol{\mathsf{y}}$ (shown as black dashed line in \Cref{fig: PCA_PCE_likelihood}). 
 \item Step two: Instead of constructing a surrogate based on the orginial output space (blue dashed line in \Cref{fig: PCA_PCE_likelihood}), the PCA technique is used to reduce the size of the high-dimensional output space. Then PCA and PCE are combined to construct a surrogate.
 \item Step three: Compute error estimates $\varepsilon_{\text{DR}}$ and $\varepsilon_{\tilde{Y}}$ in \Cref{fig: PCA-PCE}. Choose an optimal principal component number $N'$ for a surrogate construction $\tilde{\mathcal{M}}(\boldsymbol{x})$ (red solid line in \Cref{fig: PCA_PCE_likelihood}).
 \item Step four: Incorporate uncertainties into the likelihood function $\mathcal{L}(\boldsymbol{x}|\mathcal{Y})$ as shown in \Cref{eq: modified_Likelihood function} and \Cref{eq: COV_uncertainties} linked with observations $\mathcal{Y}$ and surrogate model $\tilde{\mathcal{M}}(\boldsymbol{x})$ constructed in Step three.
 \item Step five: Pass the PCA-PCE surrogate into the Bayesian inversion process and run MCMC to get the samples from the posterior $\pi(\boldsymbol{x}|\mathcal{Y})$. 
\end{itemize}
Together, these steps describe a novel approach for constructing high-dimensional output surrogates and their application in a sequential Bayesian inversion framework. It is particularly noteworthy that the approach also extends the classes of regression models that can reliably predict outputs for complex, high-dimensional problems possessing correlated outputs to include PCE, which is widely used for UQ.

\section{Application}\label{sec:numerics}\noindent
\subsection{Problem statement and details}\label{sec problemstatement}
\begin{figure}[htbp]
\centering
 \includegraphics[width=90mm]{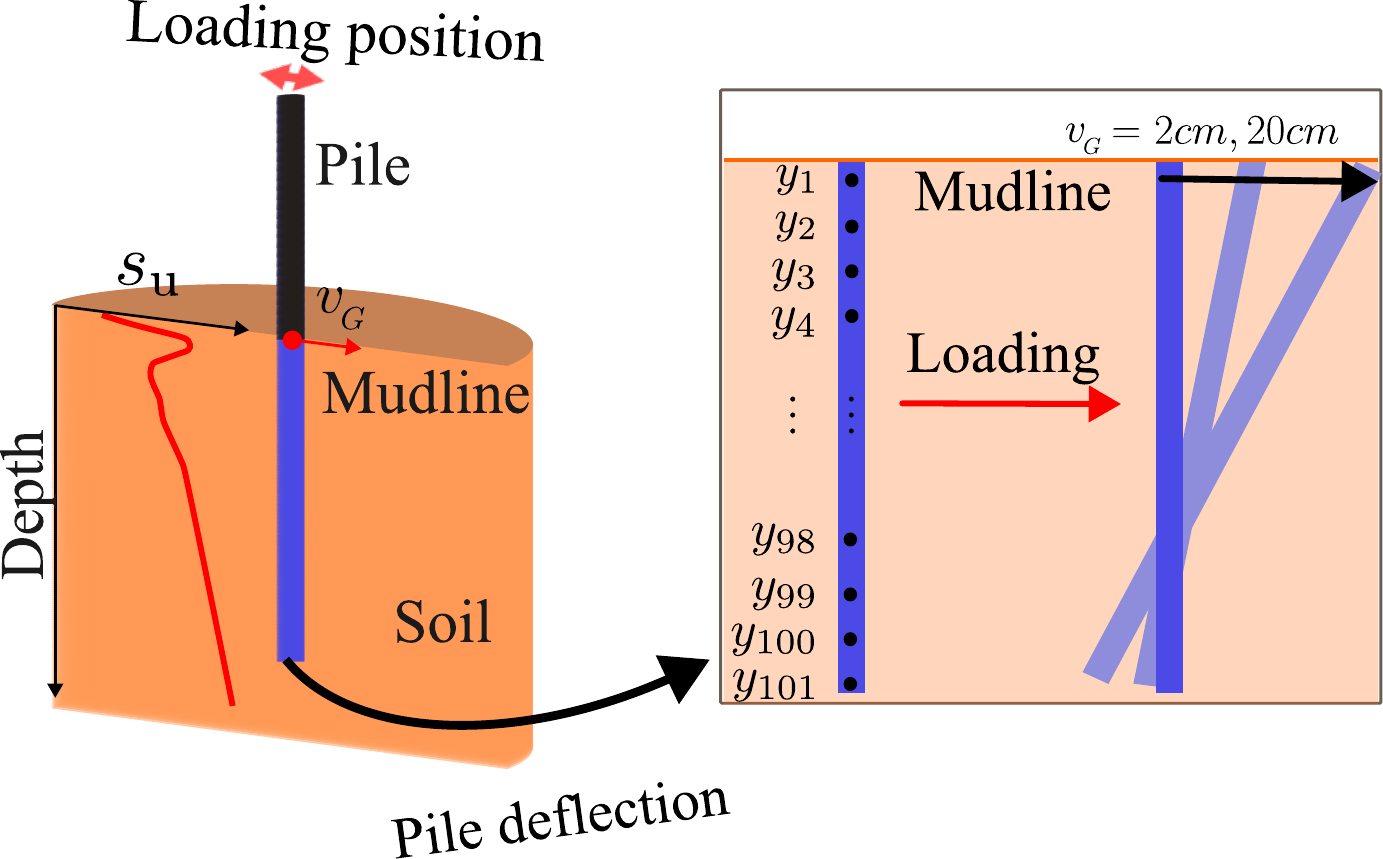}
 \caption{Model setup: Laterally loaded pile in a stiff glacial clay at Cowden}
 \label{fig: cowdensite}
\end{figure}

\noindent
To investigate the capabilities of the proposed approach, an example problem in the civil/geotechnical engineering domain, involving high-dimensional outputs, is presented. The problem is formulated based on the comprehensive field testing and numerical modelling study carried out for the PISA offshore wind turbine (OWT) project \cite{byrne2020,zdravkovic2020}. This project examined the behaviour of monopiles, supporting OWTs, when subjected to lateral loading typical of offshore environmental conditions (winds and waves). The original investigation analysed soil data from field and lab studies at the Cowden, UK glacial clay site \cite{zdravkovic2020ground} to determine parameters for the numerical model. This allowed for informed prior predictive estimates, also known as Class A (blind) predictions \cite{lambe1973}. The problem involved the 3D simulation of laterally loaded piles embedded in clay soil at the Cowden test site. All numerical analyses were conducted using the \textit{Imperial College Finite Element Program} (ICFEP) \cite{potts2001}, a deterministic forward simulator. The model's accuracy was verified by \cite{zdravkovic2020} through comparisons of predictions and field observations, including the deformed shape of the piles with depth below the ground surface, as shown in \Cref{fig: cowdensite}. 

The pile employed in this study conforms to the CL2 type outlined in \cite{zdravkovic2020}, featuring a diameter of 2.0 meters and a embedded length-to-diameter ratio of $L/D = 5.25$. It was monotonically loaded horizontally at a height of 10m above ground surface (see the schematic in \Cref{fig: cowdensite}). The finite element mesh and boundary conditions adopted are as described in \cite{zdravkovic2020}. The pile material was assumed to be linearly elastic, characterised by a Young's modulus of $E = 200$ GPa and a Poisson's ratio of $\mu = 0.30$. For the soil material, modified Cam Clay combined with a non-linear Hvorslev surface was used to characterise the glacial clay.

\begin{figure}[htbp]
\centering
 \includegraphics[width=140mm]{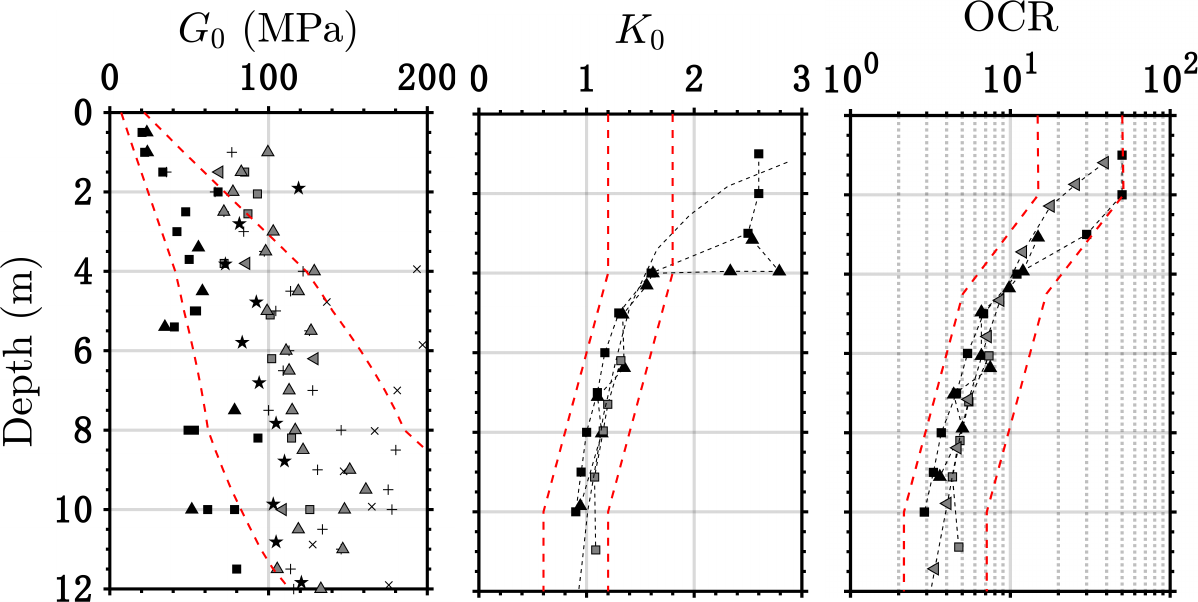}
 \caption{Parameter variations of soil properties \citep{zdravkovic2020}}
 \label{fig: cowdensite_prior}
\end{figure}

Most of the soil parameters adopted for this numerical model are presented in \Cref{table: soil-pile parameters}. The remaining parameters related to the profiles of intial elastic shear stiffness modulus, $G_0$, in the ground, the coeficient of earth pressure at rest, $K_0$, and the overconsolidation ratio, OCR. \Cref{fig: cowdensite_prior} shows the distributions of measured field and laboratory data, along with the defined upper and lower bounds of $G_0$, $K_0$ and OCR used in the PISA finite element calculations. Clearly, there is much scatter in the measured data in these plots, and consequently there is considerable uncertainty involved in selecting appropriate profiles for analysis.

In this study, the soil properties given in \Cref{table: soil-pile parameters} were assumed to be fixed, but the $G_0$, $K_0$ and OCR values were taken as uncertain and varying within the ranges shown in \Cref{fig: cowdensite_prior} (dashed red lines). Inverse analyses were performed to establish the likely values and distributions of these three parameters that gave the best estimates of the field measurements of pile deflections under lateral load. Consequently, in the forward model, these three parameters were taken as random variables sampled from the input domain $\mathcal{D}_{\boldsymbol{X}}$ through a distribution defined by the prior $\pi(\boldsymbol{x})$. The observations are taken as the measured pile deflections at 101 positions down the embedded length of the pile (each position having the same $x$ and $y$ but a different $z$ coordinate), and these are considered at two stages of the loading process, namely when the pile had a ground level displacement (at mudline) of $v_G = 2$cm and $v_G = 20$cm, respectively. The pile deflection profiles were collated in an observational data set $\mathcal{Y} = \{\boldsymbol{y}^{(i)}\}_{i=1}^2, \boldsymbol{y}^{(i)} \in \mathbb{R}^{101 \times 1}$.
\begin{table}[htbp]
\centering
\caption{Summary of model parameters for Cowden glacial clay \citep{zdravkovic2020}}
\begin{tabular}{ll}
\hline
Component & Parameters \\ \hline
Strength & $X$ = 0.548, $Y$ = 0.698, $Z$ = 0.100\\
Nonlinear Hvorslev surface-shape & $\alpha$ = 0.25, $n$ = 0.40\\
Nonlinear Hvorslev surface-plastic potential & $\beta$ = 0.2, $m$ = 1.00\\
Virgin consolidation line & $v_1$ = 2.2, $\lambda$ = 0.115\\
Nonlinear elasticity-swelling behaviour & $\kappa$ = 0.021\\
Nonlinear elasticity-shear stiffness degradation & $a$ = $\num{9.78e-5}$, $b$ = $\num{0.987}$, $R_{G,min}$ = $\num{0.05}$\\ \hline
\end{tabular}
\label{table: soil-pile parameters}
\end{table}

\begin{table}[htbp]
\centering
\caption{Priors for parameter variations}
\begin{tabular}{llll}
\hline
$\boldsymbol{x}^{\mathcal{M}}$ & $\mathcal{D}_{\boldsymbol{x}^{\mathcal{M}}}$ & Unit & Physical meaning\\ \hline
$G_{\text{0}}$ & $\mathcal{U}${(}55000,165000{)} & kPa & Maximum shear modulus \\ 
$K_{\text{0}}$ & $\mathcal{U}${(}1.2,1.8{)} & -- & Adopted variation with depth of the at-rest coefficient of earth pressure\\
OCR & $\mathcal{U}${(}15,50{)} & -- & Overconsolidation ratio \\
\hline
\end{tabular}
\label{table: input parameters}
\end{table}

Predictive estimation was carried out by quantifying the uncertainty of the three soil parameters $G_0$, $K_0$ and OCR in reproducing the observations, $\mathcal{Y}$. Profiles of the soil properties were parameterised to be a function of a single scalar value as shown in \Cref{fig: cowdensite_prior}. The training dataset was constructed assuming the parameters are independent and uniformly distributed between the bounds shown in \Cref{table: input parameters}, where the values are based on the ranges given in \Cref{fig: cowdensite_prior}. Latin hypercube sampling was employed to sample from the input space $\mathcal{D}_{\boldsymbol{x}^{\mathcal{M}}}$. The \textit{uncertainty parameters} $\boldsymbol{x^{\varepsilon}}$ were redefined at each of the two stages, following a non-informative uniform prior such that $\pi(\boldsymbol{x^{\varepsilon}})=\sigma^2 \sim \mathcal{U}(0,\text{max}(\mathcal{Y}_i))$, where $\mathcal{Y}_i = \{y^i_1,\cdots,y^i_{101}\}$.

\subsection{DR-based surrogate}\noindent
To handle the high-dimensional pile deflection outputs, following \Cref{fig: PCA-PCE}, it is preferred to extract the features and build a surrogate directly in the reduced spaces. As described above, the number of principal components (PC) to be used, denoted as $N'$, is determined based on $\varepsilon_{\text{DR}}^{threshold}$. The selection criteria for $N'$ is outlined in \Cref{fig: PCA-flowchart}, using a 2\% reconstruction error threshold and initially achieved for $N'$ = 1. This implies that the original pile deflection output $\textbf{y} \in \mathbb{R}^{101 \times 1}$ can be sufficiently represented using one principal component, $\boldsymbol{z} \in \mathbb{R}^{1 \times 1}$. The result is that each call of the likelihood function now only needs to evaluate one PCE function to carry out the MCMC simulation for this high-dimensional inverse problem, thereby improving computational efficiency by 2 orders of magnitude. The fact that only one principal value is needed in this context is due to the correlation among pile deflections along its length, which is influenced by the relative stiffness between the pile and the surrounding soil.

Based on the principal component obtained above, a PCA-PCE structure is constructed following \Cref{sec: PCA-PCE steps}. Like any other surrogates, the quality of the proposed surrogate is highly dependent on the amount of the training dataset used. Here the performance of the PCA-PCE surrogate is assessed against test runs using the Mean Absolute Percent Error (MAPE) defined as:
\begin{equation}
 \text{MAPE} = \frac{1}{K'} 
\left ( \sum_{i=1}^{K'}\frac{\mathcal{M}(\boldsymbol{x}) - \tilde{\mathcal{M}}(\boldsymbol{x})}
{\mathcal{M}(\boldsymbol{x})} \right ) \times 100
\end{equation}
where $K'$ is the number of FE runs for the test dataset, $\mathcal{M}(\boldsymbol{x})$ is the FE results and $\tilde{\mathcal{M}}(\boldsymbol{x})$ is the prediction based on PCA-PCE surrogate. 

\begin{figure}[htbp]
\centering
 \includegraphics[width = 90mm]{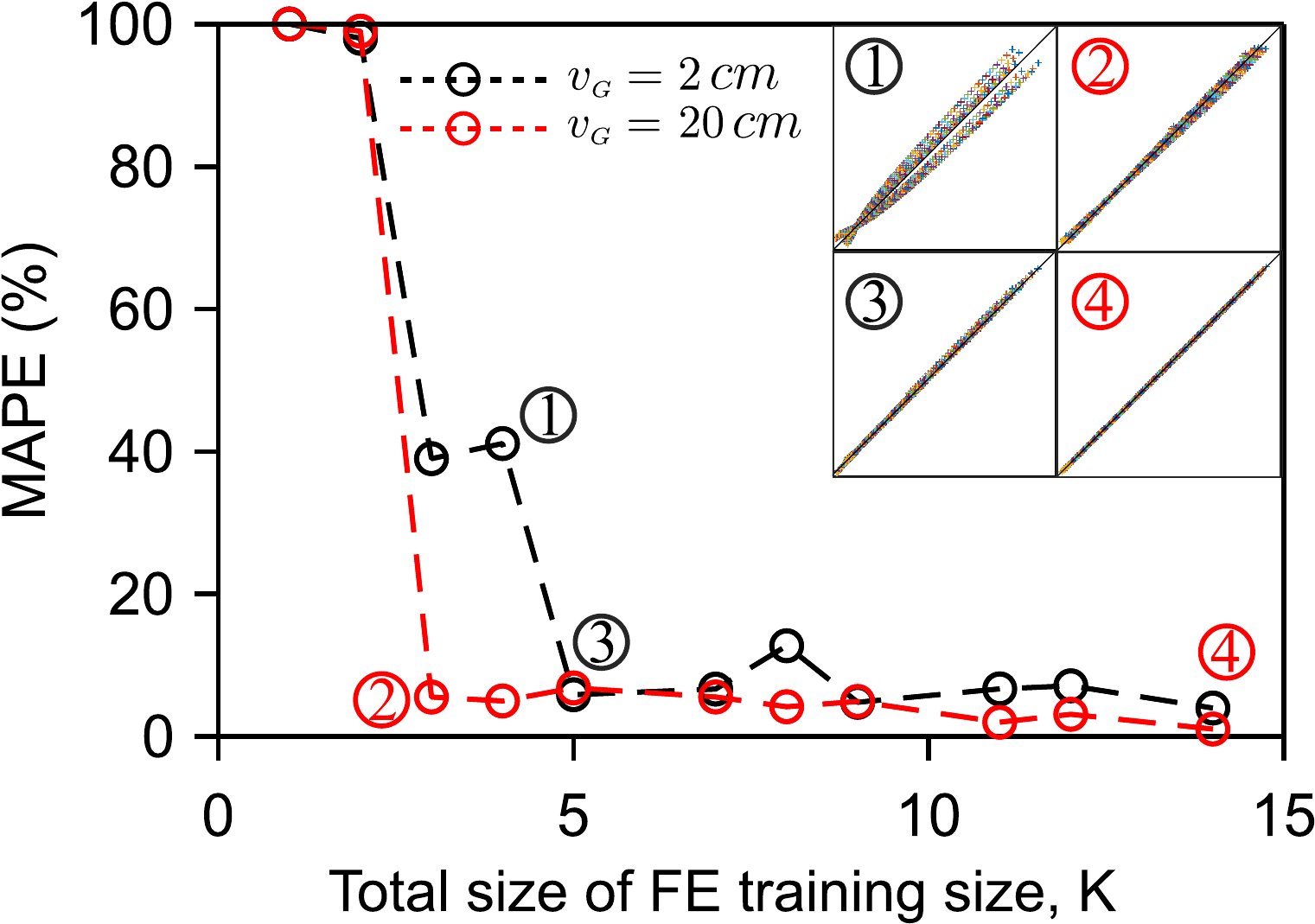}
\caption{FE training size against MAPE}
\label{fig: surrogateruns}
\end{figure}

The black and red lines in \Cref{fig: surrogateruns} represent the MAPE error for the two monopile loading stages (i.e., reaching the lateral ground movement of 2 cm and 20 cm). Only 14 forward finite element simulations were required to construct a reliable surrogate model with MAPE of less than 5\%, while 7 additional runs were performed for cross-validation testing. For both stages, the MAPE was found to decrease with an increasing number of training runs, as expected. To visualise the performance of the PCA-PCE surrogate, four cross-validation tests are presented to compare the fit between the original FE outputs and the surrogate predictions. Notably, the outputs of stage one exhibit greater non-linearity compared to those of stage two, necessitating more FE runs to achieve the same MAPE accuracy.

\subsection{Inverse analysis of a laterally loaded pile}\noindent
The effectiveness of the inverse analysis is dependent on the performance of the surrogate model within the Bayesian inversion framework. Ideally, similar to the selection of the principal components, the number of realisations in the training set should be objectively defined based on a \textit{surrogate error threshold} $\varepsilon_{\tilde{\boldsymbol{Y}}}^{threshold}$ (e.g., 5\% MAPE). However, in practice, computational constraints often make it impractical to use a predefined $\varepsilon_{\tilde{\boldsymbol{Y}}}^{threshold}$, especially when a single computational simulation can take hours to days to run. To demonstrate the value of the proposed DR-based surrogate modelling approach in scenarios where FE runs are prohibitively expensive, inverse analysis of the laterally loaded pile described in \Cref{sec problemstatement} is considered for the four cases, listed in \Cref{table: AnalysisCases}. These cases collectively explore the impact of training dataset size and the application of dimensional reduction on the surrogate MAPE. 

\begin{figure}[htbp]
\centering
 \includegraphics[width = 140mm]{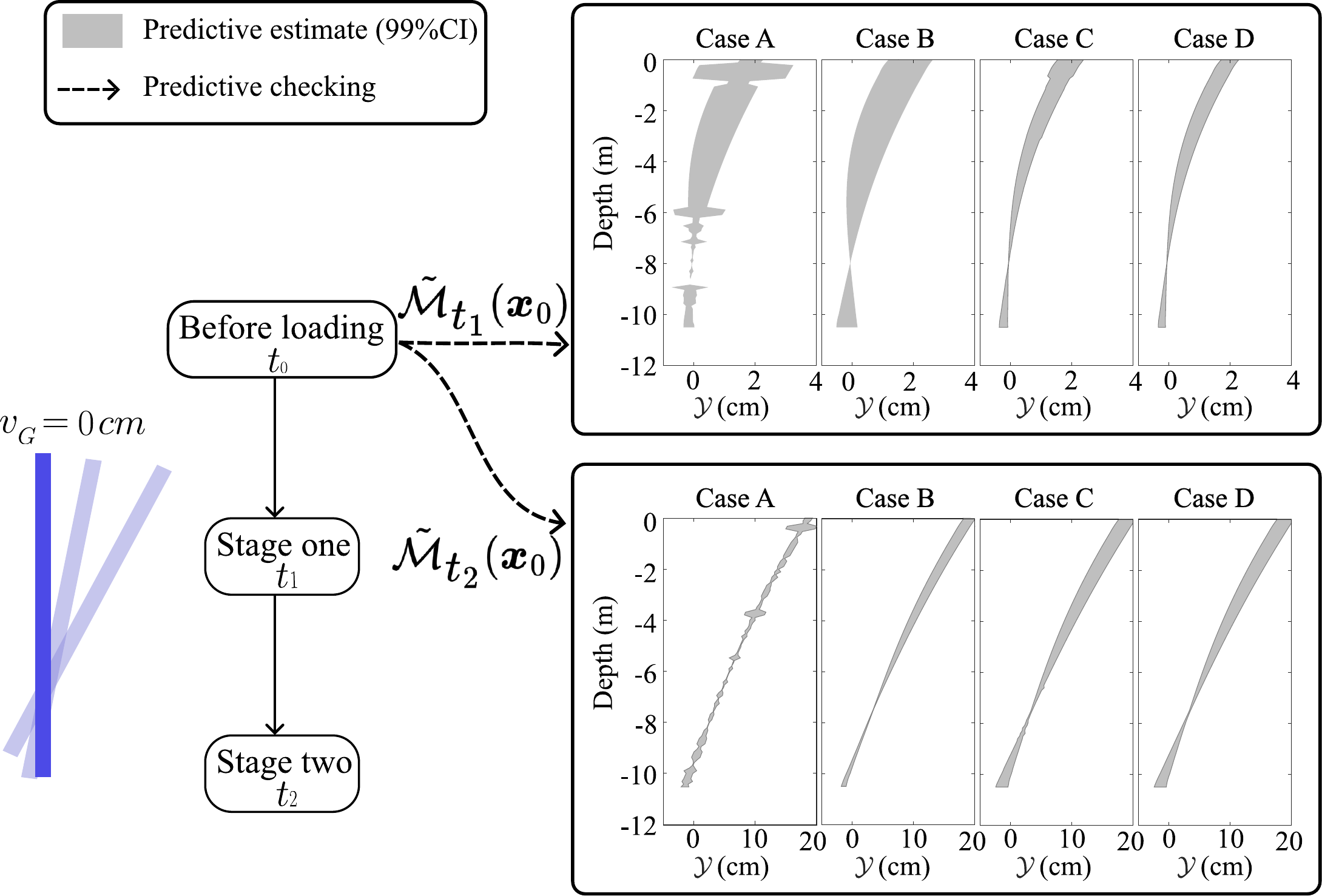}
\caption{Predictive estimates based on initial prior $\pi(\boldsymbol{x}_0)$}
\label{fig: FE_results_prior}
\end{figure}

\begin{figure}[htbp]
\centering
 \includegraphics[width = 140mm]{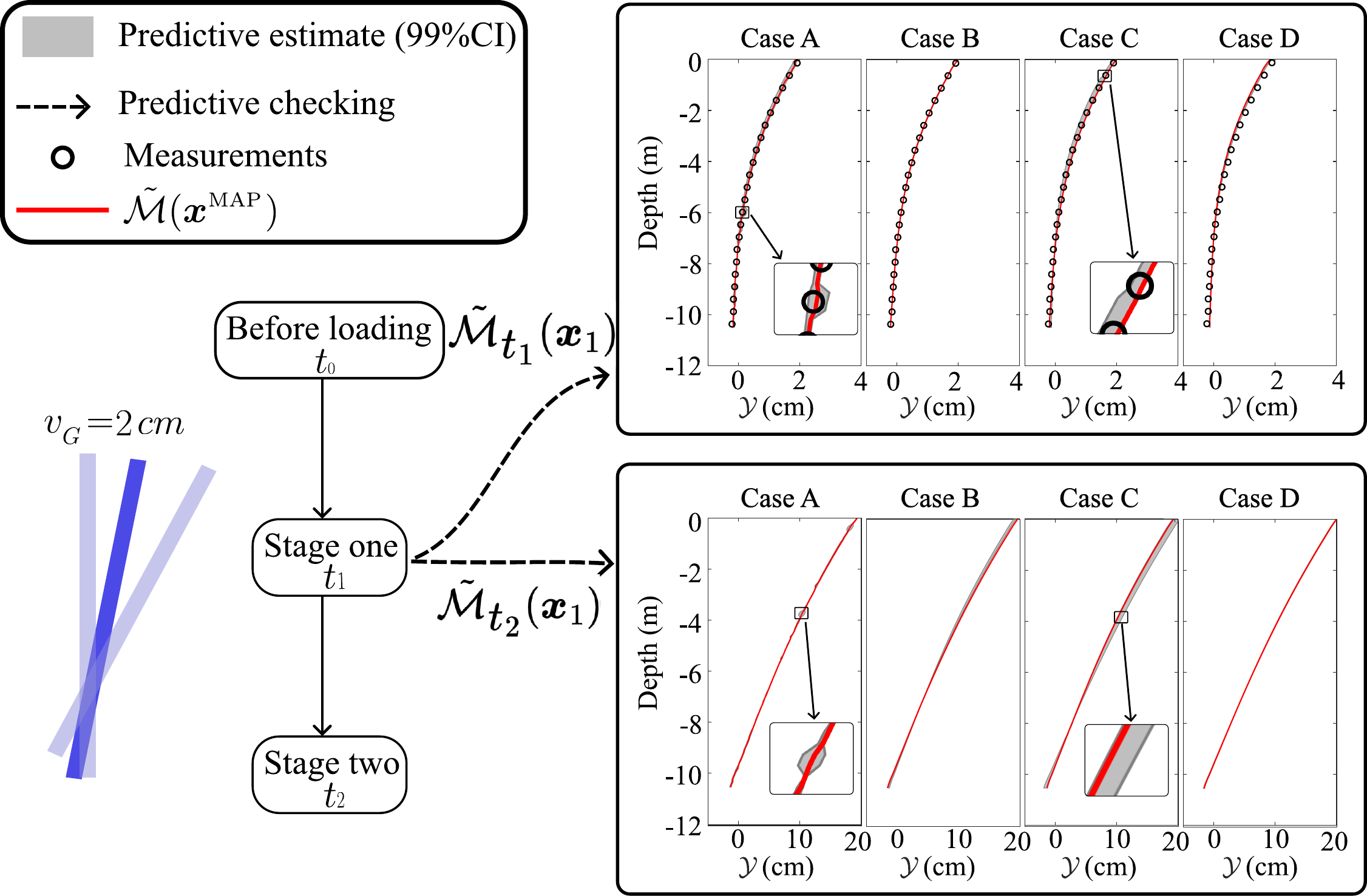}
\caption{Model checking based on $\pi(\boldsymbol{x}_1)$ (current stage $t=1$)}
\label{fig: FE_results_stageone}
\end{figure}

\begin{figure}[htbp]
\centering
 \includegraphics[width = 140mm]{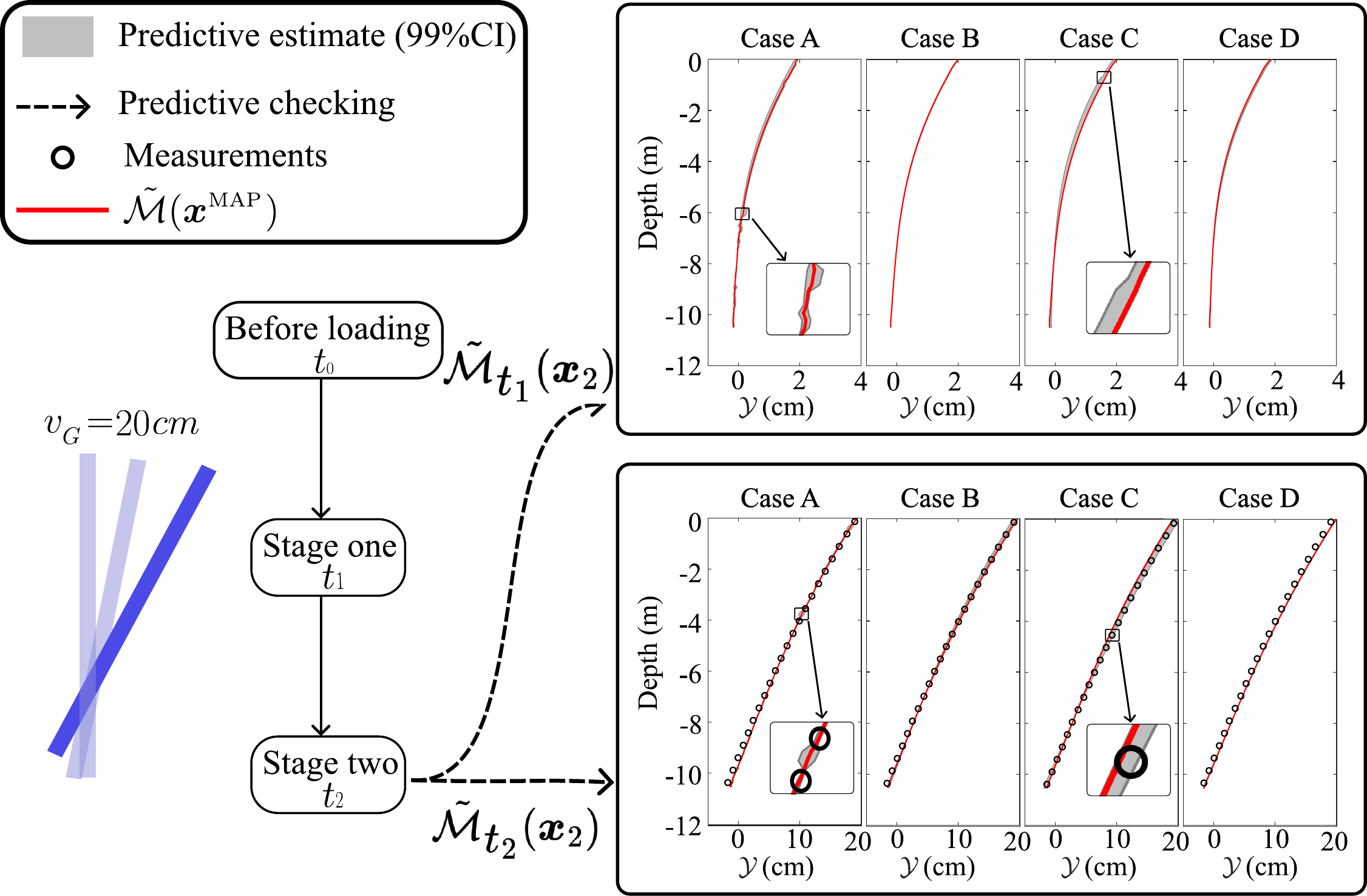}
\caption{Model checking based on $\pi(\boldsymbol{x}_2)$ (current stage $t_2$)}
\label{fig: FE_results_stagetwo}
\end{figure}

Predictive estimates, presented as 99\% confidence intervals (CI), are calculated before the start of the project (prior predictive). This stage is labelled as $t_0$ and uses the DR-based surrogate model along with Bayesian inversion to obtain forward predictions. In this sense, it is not involving inverse calculation as it does not use any of the field observations. Predictive estimates, along with the realisation of the maximum posterior, $\boldsymbol{x}^{\text{MAP}}$, are calculated using inverse analysis with the field data collection (posterior predictive) ($v_G$ = 2 cm and $v_G$ = 20 cm). These two stages are labelled as$t_1$ and $t_2$. Predictive estimates serve as a means to evaluate the spread of the surrogate model response, $\tilde{\mathcal{M}}(\boldsymbol{x})$ at a particular stage given a set of parameter distributions. The results of each of the three stages, $t_0$ to $t_2$, are presented in \Cref{fig: FE_results_prior,fig: FE_results_stageone,fig: FE_results_stagetwo}, respectively. In the figures, predictive estimates are indicated by the dashed arrow with the base of the arrow being the current stage and the tip of the arrow being the stage being estimated. Both forward and backward predictive estimates can be performed for both the forecast and the hindcast model response. 

When comparing the results from all cases, it can generally be said that irrespective to the use of a DR technique or the size of the training dataset, uncertainties in model parameters were reduced following the assimilation of observational data. However, obvious discrepancies in the quality of predictive estimates are noted when comparing cases with and without DR and using surrogates trained on different training datasets. Three key results from this investigation are reported herein. 

Firstly, with a low number of FE simulations in the training set for Cases A and B, it is unsurprising that the uncertainties in predictive estimates were large for both, particularly evident in the prior predictive. However, Case A stands out due to its notably poor performance, yielding unrealistic deflection curves at all stages. This is attributed to representing the pile deflection outputs as unconnected PCEs, each fitting an individual deflection observation (QoI) without any correlation with neighboring observations. In Case B, dimensionality reduction (DR) is used to embed the high-dimensional pile deflection outputs into a lower-dimensional subspace before constructing the PCEs. This approach trains the PCEs on a reduced space that retains non-linear features. This approach overcomes the limitations of independent sub-surrogates with limited training data, ensuring the covariance matrix constrains predictions within a realistic range and enhances predictive accuracy. Even as predictive uncertainties decreased in stages 1 and 2 with additional observations, Case A continued to show instability in providing accurate pile deflection predictions (see magnification window).

\begin{table}[H]
\caption{Analysis cases considered in this study}
\label{table: AnalysisCases}
\centering
\begin{tabular}{cccc}
\hline
Case & FE runs & Apply DR & MAPE (two stages) \\ \hline
A & 3 & \XSolidBrush & 39.8\%, 5.5\%\\
B & 3 & \Checkmark & 24.5\%, 4.8\% \\
C & 14 & \XSolidBrush & 4.5\%, 1.1\%\\
D & 14 & \Checkmark & 5.0\%, 0.9\% \\ \hline
\end{tabular}
\end{table}
Secondly, as shown in Case C in \Cref{fig: FE_results_prior,fig: FE_results_stageone,fig: FE_results_stagetwo}, despite some regions of predictive estimates remain relatively poor, using a larger training dataset helps to constrain the uncertainties in the predictive estimate. This improved performance is attributed to a greater exploration of the training dataset domain $\mathcal{D}_{\boldsymbol{X}}$ and a reduction in prediction uncertainty. With increasingly larger training datasets, the correlation between individual points in the high-dimensional outputs can be expressed as the sum of the products of the individual coefficients of the basis functions. However, compared with Case D, this result is approached asymptotically and computational costs increase with the dimensions of the outputs. Therefore, it remains advantageous to leverage the information present in the original dataset, particularly features such as the covariance between points, to construct robust surrogates with fewer realisations of the forward model.

Thirdly, results from this study show that for an equivalent training dataset of either a small or a large number of FE runs, a better high-dimensional surrogate model can always be expected when DR techniques are employed prior to the construction of the PCEs. It is worth noting, however, that the results presented in this study are sensitive to the uncertainty term defined in \Cref{eq: modified_Likelihood function}. While this study assumes that the residual matrix is diagonal, as is the case where uncertainty only represents \textit{observation errors}, it may be equally valid to assume that the uncertainty term takes the form of a non-diagonal covariance matrix. However, this term is difficult to quantify without detailed knowledge of the systematic modelling error, such cases are beyond the scope of this paper.
\section{Summary and Conclusions} \label{sec:conclusion}\noindent
Dimensionality reduction (DR) and surrogate modelling are both key ingredients for rapid, high-dimensional uncertainty quantification. This study presents a two-step strategy to combine these techniques to handle high-dimensional outputs, named DR-based surrogate modelling. Specifically, this study uses a combination of \textit{principal component analysis} and \textit{polynomial chaos expansion} to validate the performance of a high-dimensional surrogate model in predicting the response of a full-scale pile subject to lateral loading. The proposed surrogate technique is then tested for its ability to perform Bayesian inversion with and without DR techniques and for a sufficient and insufficiently trained model.

The DR-based surrogate is shown to outperform traditional PCE based surrogates, which are directly constructed for each output point (Quantity of Interest). It was shown that through dimensionality reduction, the approach is able to reduce the prevalence of unrealistic predictive estimates. This approach is shown to be particularly valuable in scenarios with inadequate training datasets, suggesting that DR-based surrogates can offer valuable insights into the structure of the high-dimensional outputs and in scenarios where the forward model is expensive to compute. Moreover, the versatility of the proposed DR-based surrogate framework in the UQ process is noteworthy, as it is not restricted to any specific problem.

In future extensions of this work, focus will be given to: (1) exploration of different combinations of surrogates and dimensionality reduction (DR) techniques; (2) extension of DR-based surrogates to the input space; (3) inclusion of various sources of uncertainties within the uncertainty quantification framework alongside the proposed surrogate.

\section{Software and data availability}\noindent
All code related to sequential Bayesian inversion and PCA-PCE surrogate is based on MATLAB software. MCMC is performed using an open-source toolbox UQLAB available at \href{https://www.uqlab.com/}{https://www.uqlab.com/}. The data underlying this paper will be shared on reasonable request to the corresponding author.

\newacronym{mathcaly}{$\mathcal{Y}$}{set of observation vectors}
\newacronym{la}{L}{Los Angeles}
\newacronym{un}{UN}{United Nations}

\glsaddall

\clearpage
\newpage

\ifnum\classstyle=0
\printbibliography
\fi
\ifnum\classstyle=1
\printbibliography
\fi
\ifnum\classstyle=2
\printbibliography
\fi
\ifnum\classstyle=3
\printbibliography
\fi

\newpage
\appendix

\ifnum\online=0
\fi
\ifnum\online=1
\fi

\end{document}